\documentclass[12pt]{article}
\usepackage{graphicx}
\usepackage{amsmath,amsthm,amsfonts,amssymb}
\input{xy}
\usepackage[all]{xy}

\newtheorem{thm}{Theorem}[section]
\newtheorem{cor}[thm]{Corollary}
\newtheorem{lem}[thm]{Lemma}
\newtheorem{pro}[thm]{Proposition}
\newtheorem{defi}[thm]{Definition}
\newtheorem{rem}[thm]{Remark}
\newtheorem{ejem}[thm]{Example}

\def\qed{\begin{flushright} \QED \end{flushright}}

\newcommand\ZZ{{\mathbb{Z}}}

\newcommand{\Ga}{\Gamma}

\def\B{{\mathcal B}}
\def\C{{\mathcal C}}
\def\D{{\mathcal D}}

\def\F{{\mathcal F}}

\def\U{{\mathcal U}}

\def\qed{\hfill \mbox{$\square$}\bigskip}

\def\ker{\mathop{\rm Ker}\nolimits}

\def\Im{\mathop{\rm Im}\nolimits}

\def\Aut{\mathop{\rm Aut}\nolimits}

\def\st{\mathsf{St}}
\def\Gal{\mathsf{Gal}}

\def\aut{\mathsf{Aut}}

\begin{document}

\sf

\title{Connected gradings and fundamental group}
\author{Claude Cibils, Mar\'\i a Julia Redondo and Andrea Solotar
\thanks{\footnotesize This work has been supported by the projects  UBACYTX212, PIP-CONICET 5099, and CONICET-CNRS.
The second and third authors are  research members of
CONICET (Argentina). The third author is a Regular Associate of ICTP Associate Scheme.}}

\date{}

\maketitle

\thanks{}
\begin{abstract}
The main purpose of this paper is to provide explicit computations of the fundamental group of several algebras.
For this purpose, given a $k$-algebra $A$, we consider the category of all connected gradings of $A$ by a group $G$
and we study the relation between gradings and Galois coverings.
This theoretical tool gives information about the fundamental group of $A$, which allows its computation using complete lists of gradings.
\end{abstract}



\section{\sf Introduction}

The main goal of this article is to provide explicit computations of the intrinsic fundamental group of some algebras.
For this, we study in detail the relation between gradings and Galois coverings of the algebra considered as a $k$-linear category with one object.
Particular attention is paid to matrix algebras, since the problem of classifying gradings of these algebras has been extensively treated in the
literature (see \cite{al,baseza,baza,boboc,bodas,bodas06,bodas07,caedasnas,chunlee,dinr,kabodas}, and also \cite{bashe}).

We recall that the intrinsic fundamental group of an algebra has been defined in \cite{cireso} using Galois coverings.
We make use of an equivalence between the category of Galois coverings and its full subcategory with objects
obtained from the smash product construction, which is deeply attached to connected gradings. We replace algebras by linear
categories over a base ring: a category over a ring $k$ is considered as an algebra with several objects, see \cite{mi}, and a
$k$-algebra $A$ can be viewed as a $k$-category with a single object and endomorphism ring equal to $A$.
Note that in \cite{green,grma} a relation between gradings and coverings is established for quivers with relations. In this paper we
consider an intrinsic context, i.e., the categories are not given by a presentation.

When performing  the computation of the fundamental group of an algebra, one faces the problem of classifying and organizing their connected gradings. The methods we introduce allow the computation of the fundamental groups of matrix algebras, triangular matrix algebras, group algebras and diagonal algebras. We restrict to connected gradings and we prove that the matrix algebras do not admit a universal grading.
Indeed, there exist at least two non-isomorphic Galois coverings or, equivalently, two non-isomorphic connected gradings which are
simply connected, in the sense that they have no nontrivial Galois coverings. In particular this provides a confirmation of the fact that the fundamental group of an algebra takes into account the matrix structure, in other words it is not a Morita invariant.

In Section 2 we show that the connectedness of gradings is the right notion which corresponds to the connectedness of the smash product associated.
We recall the concept of Galois covering and we observe that the smash product construction gives examples of Galois coverings. We
describe in detail the morphisms between smash coverings.

In Section 3 we make an explicit comparison between Galois coverings and smash coverings of a $k$-category $\B$. More precisely,
we provide an equivalence between the category $\Gal(\B, b_0)$ of Galois coverings of $\B$ and its full subcategory
$\Gal^\#(\B, b_0)$, whose objects are the smash product coverings. We consider the fundamental group that has been defined in
\cite{cireso} using Galois coverings and show that we can restrict to  smash coverings when computing the fundamental
group $\pi_1(\B,b_0)$.

In the following sections we focus on the description of connected gradings of certain algebras in order to compute their fundamental group.
As a rule, we wonder about the existence of a universal grading, since
when such a grading exists the grading group is isomorphic to the fundamental group of the algebra.

In Section 4 we consider matrix algebras: we prove that there is no universal covering by providing two non-isomorphic simply connected gradings.
Despite the fact that they appear to be very different in nature, we show that they have a unique largest common nontrivial quotient.
Using the classification of gradings of $M_2(k)$ given by C. Boboc, S. D\u asc\u alescu and R. Khazal \cite{kabodas} and of $M_3(k)$
given by C. Boboc, S. D\u asc\u alescu \cite{bodas07}, we compute the fundamental group of these algebras in case the field is
algebraically closed of characteristic  different from $2$ and $3$, respectively.
Using analogous methods and the classification of Yu. A. Bahturin and M.V. Zaicev \cite{baza},
we compute the fundamental group of $M_p(k)$, where $p$ is prime and $k$ an algebraically closed field of characteristic zero, which is the direct product of the free group on $p-1$ generators with the cyclic group of order $p$. Finally we compute the fundamental group of triangular matrix algebras, using results of A. Valenti and M.V. Zaicev in \cite{vaza},
without any hypothesis on the characteristic of the field $k$. The fundamental group in this case is the free group on $n-1$ generators.

In Section 5 we first prove that the natural grading of a group algebra is simply connected. Next we consider in detail the
group algebra of the cyclic group of order $p$, where $p$ is a prime, in case of a field of characteristic $p$.
This algebra is isomorphic to the truncated polynomial algebra $k[x]/(x^p)$, and we show that it does not admit a universal grading.
Nevertheless, we provide a complete description of its connected gradings, and we conclude that the fundamental group of the truncated
polynomial algebra in characteristic $p$ is the product of the infinite cyclic group and the cyclic group of order $p$.

Finally, in Section 6 we consider the group algebra $kG$, for $G$ an abelian group of order $n$ and $k$ a field with enough $n$-th roots
of unity or, equivalently, the algebra $k^E$ of all maps from $E$ to $k$, where $E$ is a set with $n$ elements.
In case $n$ is not square free, we show that $k^n$ has no universal covering. A quite special case occurs for $n=2$ and $k$ a field of characteristic different from $2$: there exists a universal covering. More precisely we prove that there
is only one nontrivial group providing a connected grading of the set algebra $k^2$, namely the cyclic group of order $2$, which in turn
is the fundamental group of this algebra. We end the paper by computing the fundamental group of the set algebras $k^3$ and $k^4$, using a description of all the
gradings of $k^E$ given by S. D\u asc\u alescu in \cite{da}. In case $k$ is a field containing all roots of unity of order $2$ and $3$,
we prove that $\pi_1(k^3)=C_2\times C_3$, while if $k$ contains all roots of unity of order  $3$ and $4$,
we obtain that $\pi_1(k^4)= (C_2 * C_2)\times C_6 \times C_4\times C_2.$ A detailed study of D\u asc\u alescu's classification  and the relations among the grading groups together with the techniques presented in this paragraph should lead to the computation of the fundamental group for arbitrary diagonal algebras.

We would like to thank S. D\u asc\u alescu for a valuable exchange during the preparation of this work.

\section{\sf Gradings and coverings}

Let $k$ be a commutative ring and let $\B$ be a small category such that each morphism set
${}_y\B_x$ from an object $x$ to an object $y$ is endowed with a $k$-module structure
such that composition of morphisms is $k$-bilinear. Such a category
is called a $k$-category; note that each endomorphism $k$-module ${}_x\B_x$  is a $k$-algebra and  ${}_y\B_x$  is a ${}_y\B_y$-${}_x\B_x$-bimodule.
Each $k$-algebra $A$ provides in this way a single object $k$-category.

In \cite{cima,green} it has been shown that connected group gradings and Galois coverings
are in one-to-one correspondence. We recall the definition of these categories and,
even if they are not equivalent, we make precise the relation between them.

\begin{defi}
A \textbf{grading} $X$ of a $k$-category $\B$ by a group $\Ga$ is given by a direct sum
decomposition of each $k$-module of morphisms
\[ {}_y\B_x = \bigoplus_{s \in \Ga}X^s\left({}_y\B_x\right) \]
such that $X^t\left({}_z\B_y\right)  X^s\left({}_y\B_x \right)\subset X^{ts}\left({}_z\B_x\right)$.
The \textbf{homogeneous component} of degree $s$ from $x$ to $y$ is the $k$-module $X^s\left({}_y\B_x\right)$.
\end{defi}

Next we will consider \textbf{connected} gradings in order to establish the correspondence with Galois
coverings. We use the following notation: given a morphism $f$, its source object is denoted $s(f)$
and $t(f)$  is its target object.

We will also make use of walks. For this purpose we consider the set of formal pairs $(f,\epsilon)$ as
``morphisms with sign'', where $f$ is a morphism in $\B$ and $\epsilon \in \{-1,1\}$.
We extend source and target maps to this set as follows:
\begin{center}
$s(f,1)=s(f)$, $s(f,-1)=t(f)$, $t(f,1)=t(f)$, $t(f,-1)=s(f)$.
\end{center}

\begin{defi}
Let $\B$ be a $k$-category. A non-zero \textbf{walk} in $\B$ is a sequence of non-zero morphisms
with signs $(f_n,\epsilon_n) \dots (f_1,\epsilon_1)$ such that $s(f_{i+1}, \epsilon_{i+1})= t(f_{i}, \epsilon_{i})$.
We say that this walk goes from  $s(f_1, \epsilon_1)$ to $t(f_n, \epsilon_n)$.
\end{defi}

A non-zero walk $\alpha=(f_n,\epsilon_n) \dots (f_1,\epsilon_1)$ is called \textbf{homogeneous} if each
$f_i$ is a homogeneous morphism in the graded category $\B$.
We shall denote $\mathsf{deg}f$ the degree of a homogeneous morphism $f$.
Now we define the
\textbf{degree} of a homogeneous non-zero walk $\alpha$ as follows:
\[ \mathsf{deg}\alpha= (\mathsf{deg}f_n)^{\epsilon_n}\dots (\mathsf{deg}f_1)^{\epsilon_1}. \]
As expected, a $k$-category $\B$ is called \textbf{connected} if any two objects
of $\B$ can be joined by a non-zero walk.
Moreover, a $\Ga$-grading of $\B$ is  \textbf{connected} if given any two objects in $\B$,
and any element $g\in \Ga$, they can be joined by a non-zero homogeneous walk
of degree $g$.
Of course, if a grading of a $k$-category is connected, then the underlying category is connected.
Conversely, the following easy result holds.

\begin{lem}
Let $\B$ be a connected $k$-category equipped with a $\Ga$-grading and let $x_0$ be an object of $\B$.
Assume there exist homogeneous walks of any degree from $x_0$ to itself. Then the grading is connected.\end{lem}

The definition of a connected grading restricts to algebras as follows. First recall that the  \textbf{support} of a grading $X$ of a $k$-algebra $A$ by a group $\Ga$ is:
\[ \mathsf{Supp} X=\{ s \in \Ga \ |\  X^sA \neq 0\}. \]

If the category has only one object, the following result describes the notion of connected grading of an algebra.
Note that D\u asc\u alescu \cite{da} gives the name  \textbf{faithful} to this kind of gradings.

\begin{pro}\label{2.4}
Let $A$ be $k$-algebra and $X$ be a $\Ga$-grading of $A$. The grading is connected if and only if
$\mathsf{Supp} X$ is a set of generators of $\Ga$.
\end{pro}
\noindent \textbf{Proof:} Consider the $k$-category $\B_A$ with a single object $*$ and
${}_*(\B_A)_*=A$. Assume that the grading is connected. Then, for any element $g$ of $\Ga$
there is a homogeneous non-zero walk
$\alpha=(f_n,\epsilon_n) \dots (f_1,\epsilon_1)$ such that $\mathsf{deg}\alpha=g$, which precisely means that
$\mathsf{Supp}X$ generates $\Ga$. Conversely, let $g\in \Ga$. Since $\mathsf{Supp}X$ generates $\Ga$, we
have that $g=g_n^{\epsilon_n}\dots g_1^{\epsilon_1}$ where $g_i \in \mathsf{Supp}X$
and ${\epsilon_i} = \pm 1$. Let $a_n, \dots, a_1$ be non-zero homogeneous elements
of $A$ such that $\mathsf{deg}a_i=g_i$.
Then $(a_n,\epsilon_n) \dots (a_1,\epsilon_1)$ is a non-zero closed homogeneous walk from $*$ to itself,
of degree $g$.\qed

\begin{rem}
Clearly each $\Ga$-grading of an algebra provides a unique connected grading by restricting $\Ga$ to the subgroup
generated by the support.
\end{rem}

We recall now the smash product category associated to a grading,
as defined in \cite{cima}.
This construction is compatible with the one in the algebra case, in the sense that for a finite group
$\Ga$ and  a  $\Ga$-graded algebra $A$, we recover the smash product $A\# \Ga$ given in \cite{cima}.

\begin{defi}
Let $X$ be a $\Ga$-grading of the $k$-category $\B$. The objects of the \textbf{smash product category} $\B\# \Ga$
are $\B_0 \times \Ga$ while the module of morphisms from $(b,g)$ to $(c,h)$ is
 $X^{h^{-1}g}{}_c\B_b$. In other words, morphisms are provided by homogeneous components, and composition in $\B\# \Ga$
is given by the original composition in $\B$.
The composition of morphisms is well-defined as an immediate
consequence of the definition of a graded category.

\end{defi}

\begin{rem}
Consider the previous definition for a single object $k$-category  $\B_A$ associated to a $k$-algebra $A$, and
denote $A\# \Ga = \B_A \# \Ga$. Then, the set of objects of
$A\# \Ga$ is $\Ga$, while the morphisms from $g$ to $h$ are the homogeneous elements
of degree $h^{-1}g$. If $\Ga$ is finite, the matrix algebra obtained as the direct sum of all the morphisms of this category
is precisely the smash product algebra of \cite{cima}.
\end{rem}

\begin{pro}\label{smashconnected}
 $\B \# \Ga$ is a connected category if and only if the $\Ga$-grading of $\B$ is connected.
\end{pro}
\noindent \textbf{Proof:} Note first that there is a canonical functor
$F: \B \# \Ga \to \B$ given on objects by $F(b,g)=b$, while, on morphisms,
$F$ is the inclusion map of homogeneous components.
Assume that $\B \# \Ga$ is connected and let $b$ and $c$ be objects of $\B$ and $g$ in $ \Ga$.
Consider the objects $(b,1_{\Ga})$ and $(c,g)$ in $\B \# \Ga$.
Let $\alpha=(f_n,\epsilon_n) \dots (f_1,\epsilon_1)$ be a non-zero walk from $(b,1_{\Ga})$ to $(c,g)$.
Each $f_i$ is a homogeneous morphism in $\B$, by definition of $\B \# \Ga$. \\
Note also that the target in $\B \# \Ga$ of $(f_1,\epsilon_1)$ is $(t(f_1, \epsilon_1),(\mathsf{deg}f_1)^{-\epsilon_1})$.
Moreover,
the target in $\B \# \Ga$ of $(f_2,\epsilon_2)$ is $(t(f_2, \epsilon_2),(\mathsf{deg}f_1)^{-\epsilon_1}.(\mathsf{deg}f_2)^{-\epsilon_2})$.
In this way we obtain
\[g= (\mathsf{deg}f_1)^{-\epsilon_1}(\mathsf{deg}f_2)^{-\epsilon_2}\dots (\mathsf{deg}f_n)^{-\epsilon_n},\]
consequently  $\alpha$ is a homogeneous non-zero walk from $b$ to $c$ of degree $g$.\\
Conversely, assume that the $\Ga$-grading of $\B$ is connected.
Let $(b,g)$ and $(c,h)$ be objects of $\B \# \Ga$, and consider $\alpha= (f_n,\epsilon_n) \dots (f_1,\epsilon_1)$ a homogeneous non-zero
walk in $\B$ from $b$ to $c$ of degree $h^{-1}g$.
Then $\alpha$ provides a non-zero walk from $(b,g)$ to $(c,h)$.\qed

\medskip

Coverings of $k$-categories have been introduced by K.~Bongartz and P.~Ga\-bri\-el in \cite{boga}
in order to study representation theory. We recall the definition given in \cite{cireso}.
First we define the star $\st_{b_0}\B$  at an object $b_0$ of a $k$-category $\B$
as the direct sum of all $k$-modules of morphisms with source or target $b_0$.
We note that a $k$-functor $F:\C \rightarrow \B$ induces a $k$-linear map
$F:\st_x\C \to \st_{Fx}\B$ for any object $x$ of $\C$.

\begin{defi}
Let $\C$ and $\B$ be $k$-categories. A $k$-functor
$F:\C \rightarrow \B$ is a \textbf{covering} if it is surjective on
objects and if $F$ induces $k$-isomorphisms between the
corresponding stars. More precisely for each $b_0\in\B_0$ and each
$x$ in the non-empty fibre $F^{-1}(b_0)$, the map
$$F_{b_0}^x:\st_x\C \rightarrow \st_{b_0}\B,$$
provided by $F$, is a  $k$-isomorphism.
\end{defi}

A morphism from a covering
$F:\C\rightarrow\B$ to a covering $G:\D\rightarrow\B$ is a pair of $k$-linear functors $(H,J)$
where $H: \C  \to \D$, $J: \B \to \B$ are such that $J$ is an isomorphism, $J$ is the identity on objects and $GH=JF$.

We will consider within the group of automorphisms of a covering $F:\C \rightarrow \B$, the subgroup
$\aut_1 F$ of invertible endofunctors $G$ of $\C$ such that $FG=F$.

Let $b\in \B$ and let $F^{-1}(b)$ be the corresponding fibre. This fibre is non-empty by definition, and
$\aut_1F$ acts freely on it, see \cite{le,cireso}.

\begin{defi}
A covering $F: \C \rightarrow \B$ of $k$-categories is a
\textbf{Galois covering} if $\C$ is connected and if $\aut_1F$ acts
transitively on some fibre.
\end{defi}

\begin{rem}
One can prove that for a Galois covering $F$, the group $\aut_1F$ acts transitively on any fibre, see \cite{le,cireso}.
\end{rem}

As an example of Galois coverings we have those coming from the smash product construction, that is, if $X$ is a $\Ga$-grading of the $k$-category $\B$, the functor $\B\# \Ga \to \B$, given by $(b,g) \mapsto b$ and the inclusion on morphisms, is a Galois covering with $\Ga$ as group of automorphisms.

It is useful to observe that the evident action of $\Ga$ on the smash product category $\B \# \Ga$ is given as follows.
The action on objects is given by the left action of $\Ga$ on itself. It is a free action.
Observe that for any $u\in \Ga$, a morphism from $(b,g)$ to $(c,h)$ is also a morphism from $(b,ug)$ to $(c,uh)$
since $h^{-1}g=(uh)^{-1}ug$.

We consider now Galois coverings together with a fixed object as follows. Given a $k$-category $\B$ and a fixed object $b_0$ of $\B$,
the objects of the category $\Gal (\B,b_0)$ are Galois coverings $F:\C \rightarrow \B$.
Morphisms are Galois covering morphisms $ (H,J):F_1 \rightarrow F_2$, where $H:\C_1 \to \C_2$ and $J:\B \to \B$ is an
isomorphism which is the identity on objects.

We have proved in \cite{cireso} that a morphism $(H,J)$ induces a unique group epimorphism $\lambda_H: \Aut_1F_1 \to \Aut_1F_2$
verifying $Hf=\lambda_H(f)H$, for all $f\in \Aut_1F_1$.

The following proposition describes morphisms of smash coverings in terms of the corresponding $\lambda$.

\begin{pro}
Let $b_0\in \B$ and let $F_1: \B \# G_1 \to \B$ and $F_2: \B \# G_2 \to \B$ be Galois coverings associated to connected gradings
$X_1$ and $X_2$ of $\B$ with groups $G_1$ and $G_2$.
Given a morphism of coverings $(H,J): F_1 \to F_2$ in $\Gal (\B,b_0)$, there exists a map $h:G_1 \to G_2$ such that $H(b_0,g)= (b_0,h(g))$
for all  $g\in G_1$. Moreover, $h$ is a $G_1$-morphism and $h(g)=\lambda_H(g)h(1)$, where $\lambda_H : G_1 \to G_2$ is the group morphism associated to $H$.
\end{pro}
\noindent \textbf{Proof:}
It is clear that $H(b_0,g)=(b_0,g')$ for some $g' \in G_2$ since $b_0=JF(b_0,g)=FH(b_0,g)$.
We denote $h(g)=g'$. \\
We have thus obtained that given $b_0 \in \B$, the morphism $H$ induces a map $$h:F_1^{-1}(b_0) \to F_2^{-1}(b_0).$$
Moreover, $F_i^{-1}(b_0)$ is a $G_i$-set ($i=1,2$), by indentifying $G_i$ with $\aut_1F_i$, and $\lambda_H$ makes $F_2^{-1}(b_0)$ a $G_1$-set, more
precisely, if $f\in G_1$ and $y\in F_2^{-1}(b_0)$, then $f\cdot y= \lambda_H(f)\cdot y$.
We assert  that $h$ is a morphism of $G_1$-sets. For this purpose, take $x\in F_1^{-1}(b_0)$ and $f\in G_1$, then
\[\begin{array} {rcl}
 (b_0,h(f\cdot x))& = &  H(b_0,f\cdot x)= Hf(b_0,x) \\
 & = &\lambda_H(f)H(b_0,x) = \lambda_H(f)(b_0,h(x)) \\
 & = & (b_0, \lambda_H(f)h(x)).
 \end{array} \]
Finally, $h(g)=h(g\cdot 1)=\lambda_H(g)h(1)$.
\qed


\section{\sf Fundamental group}

In \cite{cireso} the fundamental group of a connected $k$-category has been defined using Galois coverings.
Our purpose is to relate this fundamental group to connected gradings.
Let us recall the definition given in \cite{cireso}. Considering the fibre functor $$\Phi: \Gal(\B, b_0) \rightarrow \mathsf{Sets}$$
given by $\Phi(F)=F^{-1}(b_0)$, we have defined $\pi_1(\B,b_0)= \Aut \Phi$.

In order to study the fundamental group we introduce the full subcategory $\Gal^\#(\B,b_0)$:
 of $\Gal(\B,b_0)$ whose objects are the smash product Galois coverings $F: \B \# \Gamma \to \B$.

\begin{thm}
The categories $\Gal^\# (\B,b_0)$ and $\Gal(\B,b_0)$ are equivalent.
\end{thm}
\noindent \textbf{Proof:} It is immediate from \cite{cima} since any Galois covering $F:\C\rightarrow\B$ is isomorphic to the
Galois covering $\B \# \aut_1F \rightarrow \B$.
Note that the grading of $\B$ by $\aut_1F$ is not canonical, it depends on a choice of an object in each fibre.\qed

The following proposition shows that we can restrict to the subcategory $\Gal^\#(\B, b_0)$ of $\Gal(\B, b_0)$ when considering the
fundamental group $\pi_1(\B,b_0)$.

\begin{pro}
Let $F: \C \to \D$ be an equivalence of categories, $\Phi_\C : \C \to \mathsf{Sets}$, $\Phi_\D : \D \to \mathsf{Sets}$ such that
$\Phi_\D F = \Phi_{\C}$.  Then there exists an isomorphism $F^*: \Aut \Phi_\D \to \Aut \Phi_\C$.
\end{pro}
\noindent \textbf{Proof:} Recall that an element $\tau \in \Aut \Phi_\D$ is an invertible natural transformation,
that is, a family of invertible set maps $\tau_d : \Phi_\D (d) \to \Phi_\D (d)$ for every object $d$ in $\D$, which are compatible
with morphisms in $\D$. Since $F$ is a functor, it is clear that $F^*(\tau)$ defined by $F^*(\tau)_c=\tau_{F(c)}$ is an
element in $\Aut \Phi_\C$. \\
Let $\tau \in \Aut \Phi_\D$ be such that $F^*(\tau)=id$.  Since $F$ is dense, for any object $d$ in $\D$ there exists $c$ in $\C$ with an
isomorphism $\alpha: d \to F(c)$; the naturality of $\tau$ induces the commutative diagram
\[
\xymatrix{
\Phi_\D(d) \ar[d]^{\Phi_\D(\alpha)} \ar[rr]^{\tau_d}& &  \Phi_\D(d) \ar[d]^{\Phi_\D(\alpha)}\\
\Phi_\D(F(c)) \ar[rr]^{\tau_{F(c)}}& &  \Phi_\D(F(c)).}
\]
Since $\tau_{F(c)}=id$ for all $c\in \C$, this implies that $\tau_d=id$, and hence $\tau=id$. \\
In order to prove that $F^*$ is surjective, let $\sigma \in \Aut \Phi_\C$ and consider $\hat \sigma$ defined in the following way.
For any object $d$ in $\D$, we choose $c$ and an isomorphism $\alpha: d \to F(c)$; in case $d=F(c)$, we choose $\alpha=id$.
Now we define $\hat \sigma_d$ such that the following diagram is commutative
\[
\xymatrix{
\Phi_\D(d) \ar[d]^{\Phi_\D(\alpha)} \ar[rr]^{\hat \sigma_d}& &  \Phi_\D(d) \ar[d]^{\Phi_\D(\alpha)}\\
\Phi_\C(c)=\Phi_\D(F(c)) \ar[rr]^{\sigma_c}& &  \Phi_\D(F(c))=\Phi_\C(c).}
\]
Since $F$ is full we have that $\hat \sigma$ is a natural transformation and $F^*(\hat \sigma)=\sigma$. \qed

\begin{cor} Let $\Phi^\#: \Gal^\#(\B, b_0) \rightarrow \mathsf{Sets}$
be the functor given by $$\Phi(F: \B\# G \to \B)=F^{-1}(b_0)=G.$$  Then $\pi_1(\B,b_0)\cong \Aut \Phi^\#$.
\end{cor}

\begin{cor}
If $\B$ only admits the trivial connected grading, then $\pi_1(\B,b_0)=1$.
\end{cor}

An advantage of considering $\Gal^\#(\B, b_0)$ instead of $\Gal(\B, b_0)$ is explained by the following proposition,
which describes the automorphisms of the fibre functor.

\begin{pro}\label{sigma}
Let $\sigma \in \Aut \Phi^{\#}$ and let $G$ be a group grading the category $\B$ in a connected way.
The map $\sigma_G: G\to G$ is given by $\sigma_G(x)=xg$, where $g\in G$ is uniquely determined.
\end{pro}
\noindent \textbf{Proof:}
Consider a covering $F: \B \# G \to \B$. Each $g\in G$ induces an automorphism of the covering $F$, which is the identity on $\B$
and it is the left action  of $G$ on itself. We shall denote it $l_g$.
Given $\sigma \in \Aut \Phi^{\#}$ we get a map $\sigma_G: G\to G$. It must make the following diagram commutative:
\[
\xymatrix{
G
\ar[r]^{\sigma_G}
\ar[d]_{\widetilde{l_g}}
&
G
\ar[d]^{\widetilde{l_g}}
\\
G
\ar[r]^{\sigma_G}
&
G
}
\]
where $\widetilde{l_g}$ is induced by $l_g$.
So, for all $x\in G$, we get $\sigma_G(g_0x)=g_0\sigma_G(x)$. Taking $x=1$ we obtain  $\sigma_G(g_0)=g_0\sigma_G(1)$.
Note that $g_0$ is an arbitrary element of $G$. \qed

\section{\sf {Fundamental group of matrix and triangular algebras}}

Let $k$ be a field containing a primitive $n$-th root of unity $q$ and let $M_n(k)$
be the $k$-algebra of $n\times n$ matrices. The problem of classifying all the gradings of $M_n(k)$ is not solved.
Lists of gradings have been described by several authors \cite{al,baseza,baza,boboc,bodas,bodas06,caedasnas,dinr},
and the complete lists for $n=2$ and $n=3$ are obtained in \cite{kabodas} and \cite{bodas07}.

We will consider connected gradings of the algebra $M_n(k)$. In case of a
non-connected grading we shall restrict to the subgroup
generated by the support, in order to study the unique associated
connected grading.

We briefly recall the definition of the universal covering of a $k$-category and Theorem 4.6 of \cite{cireso}.

\begin{defi}
A \textbf{universal covering} $U:\U \to \B$ is an object in $\Gal(\B)$ such that for any Galois covering $F:\C \to \B$,
and for any $u_0\in \U_0$, $c_0\in \C_0$ with $U(u_0)=F(c_0)$, there exists a unique morphism $(H,1)$ from $U$ to $F$
verifying $H(u_0)=c_0$.
\end{defi}

\begin{thm}
Suppose that a connected $k$-category $\B$ admits a universal covering $U$.
Then $$\pi_1(\B,b_0) \simeq \aut_1 U.$$
\end{thm}

\begin{defi}
A connected $k$-category is  \textbf{simply connected} if its only connected grading is the trivial grading.
A connected grading is \textbf{simply connected} if the corresponding Galois covering is simply connected.
\end{defi}

We will prove that there is no universal cover for $M_n(k)$. Indeed there exist at least two
non-isomorphic connected gradings of $M_n(k)$ which provide simply connected Galois coverings. Recall that a covering is simply connected if it admits no proper Galois covering.

\begin{pro}\label{cncn} \cite{baseza,chunlee}
There exists a connected $C_n\times C_n$-grading of $M_n(k)$.
\end{pro}
\noindent \textbf{Proof:} The algebra $M_n(k)$ has a well-known presentation as follows:
\[ M_n(k)=k\{x,y\}/\langle x^n=1,y^n=1,yx=qxy \rangle  \]
where \[ x= \left(
\begin{array}{rrrrrr}
0 & 0 &  \cdots & 0 & 1\\
1& 0  & \cdots &0 & 0\\
0&1& \cdots & 0 & 0\\
\vdots &   & \ddots & & \vdots \\
0&0&  \cdots & 1 & 0
\end{array}
\right), \qquad
y=  \left(
\begin{array}{rrrrr}
q&0&0& \cdots &0\\
0&q^2&0 &\cdots&0\\
0&0&q^3& \cdots & 0\\
\vdots &  &  & \ddots & \vdots \\
0&0&0& \cdots & q^n\\
\end{array}
\right) , \]
with $q$ a primitive $n$-th root of unity.
We provide a connected grading of $k\{x,y\}$ by assigning degree $(t,1)$ to $x$ and
degree $(1,t)$ to $y$, where $t$ is a generator of $C_n$.
The group is abelian and the order of the generators is $n$, hence the ideal of relations is homogeneous.
Since the support coincides with $C_n\times C_n$, the grading is connected.\qed

\begin{pro}\label{q}
Let $\C$ be a $k$-category with a finite set of objects and one-dimensional vector spaces of morphisms
between any pair of objects $b$ and $c$, denoted by ${}_c \C_b=k\ {}_cf_b$,
verifying $$({}_df_c)({}_cf_b)=q_{d,c,b}\ ({}_df_b)$$
for any triple of objects of $\C$, where $q_{d,c,b} \in k^*$ are the structure constants.
Then $\C$ is simply connected.
\end{pro}
\noindent \textbf{Proof:} Let $G$ be a group providing a grading of the category  $\C$. As we  have already written,
we consider connected gradings. Since all the $k$-vector spaces of morphisms
are one-dimensional, they are homogeneous.
Let ${}_cs_b$ be the degree of ${}_c\C_b$. Note that for each object $b$
we have ${}_b\C_b=k$, hence ${}_bs_b=1$ and ${}_bs_c= {}_cs_b^{-1}$.
We assert that any non-zero homogeneous closed walk has degree $1$.
Indeed, since composition of non-zero morphisms is non-zero in $\C$, and since
${}_bs_c= {}_cs_b^{-1}$, a non-zero homogeneous closed walk at $b$ can be replaced by a non-zero endomorphism
of $b$, having the same degree.
Since endomorphisms of $b$ have degree $1$, the assertion is proved.
Recall that a grading is connected if for any pair of objects, any group element appears as the degree of
a non-zero homogeneous walk between them. Since the grading is connected, the group is trivial.\qed

{\begin{cor}
Let $\C$ be a category as above. Then $\pi_1 (\C) = 1$.
\end{cor}
Let ${}_jE_i$  be the matrix whose entries are zero, except the $(j,i)$ entry which equals $1$.
We recall that a \textbf{good grading} of a matrix algebra is a grading where the elementary matrices  ${}_jE_i$, also called \textbf{matrix units}, are homogeneous, see for instance \cite{dinr}.
Note that the $k$-category {${\cal M}_n(k)$} associated to a matrix algebra with respect to the idempotent elementary matrices
${}_iE_i$ is precisely a category as in the above proposition, where all the {structure} constants equal $1$.

Clearly good gradings of $M_n(k)$ and gradings of the $k$-category ${\cal M}_n(k)$ coincide.

\begin{cor}\label{goodgrading}
Let $G$ be a group providing a good grading of a matrix algebra, and assume that the corresponding grading of the $k$-category ${\cal M}_n(k)$ is connected. Then $G$ is trivial.
\end{cor}

\begin{rem}
A good grading by a nontrivial group $G$ of a matrix algebra $M_n(k)$ can be connected when $M_n(k)$ is
viewed as a category with a single object. This means that the support of the grading generates $G$.
Corollary \ref{goodgrading} makes precise that the
corresponding grading of the $k$-category ${\cal M}_n(k)$ will not be connected.
\end{rem}
}

\begin{thm}
The connected grading of the matrix algebra $M_n(k)$ by the group $C_n\times C_n$ of Proposition \ref{cncn} is simply connected.
\end{thm}
\noindent \textbf{Proof:}
We will prove that the Galois covering $\C= M_n(k) \# (C_n\times C_n)$ is simply connected.
The category $\C= M_n(k) \# (C_n\times C_n)$ has set of objects
$C_n\times C_n=\{ a^ib^j\ |\    0\le i,j\le n-1 \}$ and
\[ {}_{a^sb^l}\C_{a^ib^j}= X^{a^{i-s}b^{j-l}}M_n(k)=k(x^{i-s}y^{j-l}). \]
Hence the $k$-vector spaces of morphisms are one-dimensional with basis elements
$${}_{(s,l)}f_{(i,j)}=x^{i-s}y^{j-l}$$ and
\[ {}_{(u,v)}f_{(s,l)}{}_{(s,l)}f_{(i,j)}= x^{s-u}y^{l-v}x^{i-s}y^{j-l}= q^{l+i-v-s}{}_{(u,v)}f_{(i,j)} . \]
Finally, Proposition \ref{q} asserts that such categories are simply connected.\qed

\medskip

Each time a universal covering exists, the fundamental group is isomorphic to its Galois group.
Clearly, if there exist at least two non-isomorphic simply connected coverings, there is no
universal covering.

We will now show that this is the case for $M_n(k)$, i.e., there exists at least another simply connected grading of $M_n(k)$.
For this purpose we first provide another presentation of the matrix algebra, as a quotient of a path algebra.

\begin{pro}\label{presentamatrices}
Let $Q$ be the quiver with $n$ vertices labelled $1,\dots,n$, and arrows $x_i$ from $i$ to $i+1$ as well as reverse
arrows $y_i$ from $i+1$ to $i$ for for $1\leq i <n$. We denote $e_1,\dots,e_n$ the idempotents of the path algebra $kQ$
corresponding to the vertices. Let $I$ be the two-sided ideal of $kQ$ generated by $y_ix_i-e_i$ and $x_iy_i-e_{i+1}$ for $1\leq i <n$.
Then $kQ/I$ is isomorphic to $M_n(k)$.
\end{pro}
\noindent \textbf{Proof:} Consider the morphism of algebras $\varphi : kQ \rightarrow M_n(k)$ given by $\varphi(e_i)={}_iE_i$, $\varphi(x_i)={}_{i+1}E_i$ and $\varphi(y_i)={}_iE_{i+1}$, which is well defined by the universal property of path algebras, which are in fact
tensor algebras over the semisimple commutative algebra given by the length zero paths. This map is surjective since the matrices ${}_jE_i$
are clearly images of paths of $Q$. Moreover $I\subset\ker\varphi$ and $\dim_kkQ/I\leq n^2$.
\qed

Let $F_{n-1}$ be the free group on $n-1$ generators $s_1, \dots ,s_{n-1}$. First we introduce an $F_{n-1}$ grading of $kQ$ as follows:
for $1\le i \le n$, let $\deg e_i=1$, while for $1\le i \le n-1$ we set $\deg x_i = s_i$ and $\deg y_i = (s_i)^{-1}$.
The path algebra is a free algebra on the set of arrows with respect to
the semisimple subalgebra of vertices, so this provides a well defined grading of $kQ$. More precisely the degree of any path is the
corresponding product of the degrees of the arrows. Since the ideal $I$ is homogeneous with respect to this grading, we obtain a
grading of $kQ/I$, hence of $M_n(k)$. Note that this grading, considered as a grading of the algebra $M_n(k)$, that is, as a
grading of the single object category with endomorphism algebra $M_n(k)$, is connected since the generators of the free group are
in the support.

\begin{pro}\label{free}
The above $F_{n-1}$-grading of $M_n(k)$ is simply connected.
\end{pro}
\noindent \textbf{Proof:} The set of objects of $M_n(k) \# F_{n-1}$ is $F_{n-1}$. For $j>i$, let ${}_js_i=s_{j-1} \dots s_{i+1}s_i$.
There is a one-dimensional vector space of morphisms from a word $w$ in $F_{n-1}$ considered as an object of $M_n(k)\# F_{n-1}$
to each object ${}_js_iw$ with basis vector denoted ${{}_jE_i}^w$.
Similarly for $j<i$ there is a one-dimensional vector space of morphisms from $w$ to  ${}_js_i^{-1}w$ , with basis  ${{}_jE_i}^w$.
From $w$ to $w$, the $n$-dimensional vector space of morphisms has basis  $\{ {{}_1E_1}^w, \dots, {{}_nE_n}^w \} $.
Note that the endomorphism algebra of each object is the $n$-dimensional diagonal algebra
$k({{}_1\!{E_1}}^w)\times \cdots \times k({{}_nE_n}^w)$.
Consider now a grading of this category by a group $G$. Since the spaces of morphisms between different
objects are one-dimensional, they are homogeneous.
This fact implies that for each object $w$ the subvector space $k\left({{}_iE_i}^w\right)$ is homogeneous since
\[ {{}_iE_j}^{\left({}_j\!s_iw\right)}{{}_jE_i}^w = {{}_iE_i}^w. \]
Observe that an idempotent homogeneous element has necessarily degree $1$, hence each endomorphism
algebra has trivial grading (all elements have degree $1$).
As a consequence,
$$\mathsf{deg}\left({{}_iE_j}^{({}_j\!s_iw)}\right) = {\left(\mathsf{deg}{}_jE_i^w\right)}^{-1}.$$
Moreover, for $j>i$
\[ \mathsf{deg}\left({{}_jE_i}^w\right)= \mathsf{deg}\left({{}_jE_{j-1}}^{\left({}_{j-1}\!s_iw\right)}\right) \ \cdots\  \mathsf{deg}\left({{}_{i+2}E_{i+1}}^{s_iw}\right)\ \mathsf{deg}\left({{}_{i+1}E_i}^w\right) . \]
For $j<i$, the statement is analogous considering the inverses of the above degrees.
This complete description of any possible grading of the smash category shows that any closed homogeneous non-zero
walk has degree $1$. Consequently the grading is connected only if the group is trivial.\qed

The complete list of good gradings of a matrix algebra is obtained  in \cite{caedasnas}.
In order to compute the fundamental groups of matrix algebras, we make this classification explicit using Proposition \ref{presentamatrices}.

\begin{thm}
There is a one-to-one correspondence between good connected $G$-gradings of $M_n(k)$ and maps
$\{1, \dots, n-1 \} \to G$ such that the image generates $G$.
\end{thm}
\noindent \textbf{Proof:} Let $m$ be a map from $\{1, \dots, n-1 \}$ to $G$.
We obtain a grading of the algebra $kQ$ defined in Proposition \ref{presentamatrices} as before, namely $\deg({}_{i+1}E_i) = m(i)$ and $\deg({}_{i}E_{i+1})=m(i)^{-1}$. The ideal of relations of Proposition \ref{presentamatrices} is homogeneous  and we obtain a good grading
of $M_n(k)$. If the image of $m$ generates $G$, then the grading is connected.
Conversely, consider a good connected grading of $M_n(k)$ by a group $G$. The image of the map
$m:\{1, \dots, n-1 \} \to  G$ given by $m(i) = \deg({}_{i+1}E_i)$ generates $G$.\qed

Note that relaxing the connectedness  requirement for good gradings is equivalent to deleting the condition that the image of
each map $m$ generates $G$.
In \cite{caedasnas} the algebra $M_n(k)$ is viewed as the endomorphism algebra of a vector space $V$, and good gradings
are obtained from a grading of $V$, considering graded endomorphisms as homogeneous components.

\begin{defi}
The \textbf{quotient} of a $G$-grading $X$ of a category  $\B$ by a normal subgroup $N$ of $G$ is a $G/N$-grading ${X}/N$ of $\B$, where
the homogeneous component of degree $\alpha$ is
$$({X}/N)^\alpha {}_c\B_b = \bigoplus_{g\in\alpha}X^g {}_c\B_b.$$
\end{defi}

Observe that if $X$ is connected then $X/N$ is also connected.

The corresponding functor between the smash product coverings is precisely the canonical projection obtained through the quotient
of $\B\#G\rightarrow \B$ by $N$.

\begin{pro}\label{conocidas}
Any good connected $G$-grading of $M_n(k)$ is a quotient of the
$F_{n-1}$-grading considered before.
\end{pro}
\noindent \textbf{Proof:} Let $m_0 : \{1, \dots, n-1 \} \to  F_{n-1} $ be the map corresponding to this grading,
given by $m_0(i) = s_i$ and let $m : \{1, \dots, n-1 \} \to  G$ be another map such that the image of $m$ generates $G$.
Then the group homomorphism given by $s_i\mapsto m(i)$ is a surjective group morphism.\qed

We recall that a simply connected grading is a grading which is maximal in the sense that it is not isomorphic to a proper quotient of a connected grading.

\begin{pro}\label{cocientecomun}
Let $k$ be a field containing a primitive $n$-th root of unity. The grading by $C_n\times C_n$ of Proposition \ref{cncn} and the
grading by the free group of Proposition \ref{free} have a unique maximal common quotient $C_n$-grading.
\end{pro}

\noindent \textbf{Proof:} We denote by $X$ the grading by $C_n\times C_n$ and we observe that the vector space $X^1$ of homogeneous elements
of trivial degree is one-dimensional. Let $Y$ be the grading by $F_{n-1}$: observe that $Y^1$ is the $n$-dimensional subalgebra
of diagonal matrices.
Assume that $Z$ is a common quotient of $X$ and $Y$ and let $N$ be the normal subgroup of $C_n\times C_n$ which provides $Z$ as a quotient of $X$.
As $Z$ is a quotient of $Y$, clearly $Z^1$ contains at least $Y^1$, the diagonal matrices. Observe that the elementary diagonal matrices are
homogeneous for $X$, consequently their degrees must be elements of $N$ in order to become trivial.  Note that the set of
degrees of the diagonal matrices for $X$ is precisely $1\times C_n$. Hence $1\times C_n$ is the smallest subgroup of $C_n\times C_n$
which has a chance to meet a quotient of $Y$; let $N=1\times C_n$. In fact we assert that $X/N$ is already a good grading, in other words, elementary matrices are homogeneous. Indeed, consider the $n$-dimensional subvector space $E$ of $M_n(k)$ with basis $\{{}_2E_1, {}_3E_2, \dots, {}_{n}E_{n-1}, {}_1E_n\}$.
Recall that $x$ is the circulant matrix, which is the sum of all the previous basis vectors of $E$, while $y$ is the diagonal matrix made
with powers of the primitive root of unity $q$. Then the set $\{x,xy, xy^2,\dots, xy^{n-1} \}$
is clearly contained in $E$. Moreover the elements $xy^i$, for $0\le i \le n-1$, are homogeneous for the grading $X$, of different degrees $(t,t^i)$ where $t$ is the
generator of $C_n$. Hence they are linearly independent and they form a basis of $E$. Finally we observe that for $X/N$ all these elements
have the same degree $\overline{(t,1)}$, hence  $E$ is contained in the the set of homogeneous elements of degree $\overline{(t,1)}$ of  $X/N$.\\
Consequently each elementary matrix is homogeneous for $X/N$. Considering $Y$, we obtain $X/N$ as the quotient $Y/M$, where $M$
is the smallest normal subgroup of $F_{n-1}$ such that in $F_{n-1}/M$ all the generators of $F_{n-1}$ are equal, and this element
is of order $n$.\qed

\begin{thm} \label{23}
Let $k$ be an algebraically closed field. \begin{itemize}
\item [(1)] If $char(k)\neq 2$, then $\pi_1 M_2(k) \simeq \mathbb{Z}\times C_2.$
\item [(2)] If $char(k)\neq 3$, then $\pi_1 M_3(k) \simeq F_{2}\times C_3.$
\end{itemize}
\end{thm}
\noindent \textbf{Proof:}
Under these assumptions, the classifications of \cite{kabodas,bodas} show that all
gradings are good gradings or quotients of the one given by Proposition \ref{cncn}.
The latter and the grading by the free group have a common quotient described in Proposition \ref{cocientecomun}.
Recall that we have proved that all good gradings of $M_n(k)$ are quotients of the $F_{n-1}$-grading. We  now prove the first assertion.
We construct two inverse group morphisms between $\pi_1 M_2(k)$ and $F_1\times C_2$.
Let $\sigma \in \Aut \Phi^{\#}$, where $\Phi^{\#}:\Gal^{\#}(M_2(k)) \to Sets$ is the fibre functor.
Consider the good $F_1$-grading of $M_2(k)$. By Proposition \ref{sigma} the map $\sigma_{F_1}$ verifies $\sigma_{F_1}(x)=xg$
for some uniquely determined $g\in F_1$.
Analogously, the $C_2\times C_2$-grading provides  $\sigma_{C_2\times C_2}$  and an element $(t^a,t^b)$, where
$C_2=\langle t \rangle$.
The compatibility condition obtained when considering the maximal common quotient $C_2$ says that $t^a$ equals
the class of $g$ in $C_2$. We associate the pair $(g,t^b)$ to $\sigma$. \\
Conversely, given $(g,t^b)\in F_1\times C_2$, we will construct $\sigma \in \Aut \Phi^{\#}$ associated to it.
One needs to have maps $\sigma_G: G\to G$ for each group $G$ providing a connected grading of $M_2(k)$.
Using the classification of the gradings given in \cite{kabodas} and Proposition \ref{cocientecomun},
it is sufficient to describe $\sigma_{F_1}$ and  $\sigma_{C_2\times C_2}$.
Fix $\sigma_{F_1}(x)=xg$ and  $\sigma_{C_2\times C_2}(x)=x(\overline{g},t^b)$.
Note that these maps satisfy the compatibility condition.
Of course, for the other quotients $G$ of $F_1$ or of $C_2\times C_2$, the map $\sigma_G$ is uniquely determined
thanks to  the quotient compatibility conditions. \\
The proof of the second statement is completely analogous.\qed

Next we will prove a generalization of the preceding theorem for matrices of prime size. We will use the main result obtained by Y.
Bahturin and M. Zaicev in \cite{baza}. Considering an algebraically closed field $k$ of characteristic zero, Theorem 5.1 of \cite{baza}
states that any grading of $M_n(k)$ by a group $G$ is a tensor product of gradings, in the sense that there exists a
decomposition $n=n_1n_2$, a fine grading of $M_{n_1}(k)$ by a subgroup $G_1$ of order $n_1^2$ and a good $G$-grading of
$M_{n_2}(k)$ such that $M_n(k)$ is isomorphic as $G$-graded algebra to the tensor product algebra $M_{n_1}(k)\otimes M_{n_2}(k)$
which is obtained as an induced grading. Observe that the construction of an induced grading resembles to a tensor construction,
however it is well defined only in case one of the graded algebras involved is a matrix algebra with a good grading, see \cite{baza}.

\begin{pro}
Let $p$ be a prime and $k$ be an algebraically closed field of characteristic zero. Let $X$ be a maximal connected grading
by a group $G$ of $M_p(k)$.
Then either the group $G$ is isomorphic to $C_p\times C_p$ and the grading is fine as in Proposition \ref{cncn}, or the grading is a good
grading given by $m :\{1,\dots, p-1\}\to G$ such that $\Im m$ generates $G$.
\end{pro}

\noindent \textbf{Proof:} Since $p$ is a prime, Theorem 5.1 of \cite{baza} cited above shows that the grading is either good, or
fine with group of order $p^2$. We already know that good connected gradings are as described in Proposition \ref{conocidas}.
If the grading is fine, the order $p^2$ of the group is precisely the dimension of the matrix algebra, hence $\mathsf{Supp} X = G$.
Moreover, for fine gradings of matrix algebras, homogeneous non-zero elements are invertible by Corollary 2.7 of \cite{baza}.
Then we assert the group is not cyclic: indeed, if $G$ has a generator $t$ of order $p^2$, let $x$ be a non-zero element of
degree $p^2$, hence invertible. Note that $X^1 M_p(k) = k$. Hence $x^{p^2}\in k$ and $x^{p^2}\neq 0$, we can normalize
$x$ dividing it by a scalar in order to obtain $x'\in X^t M_p(k)$ such that $x'^{p^2}=1$. Then $M_p(k)$ would be
isomorphic to the group algebra of the cyclic group of order $p^2$, which is false since for instance the former is commutative. \\
Consequently a fine connected grading of $M_p(k)$ is given by $C_p\times C_p$. As before $\mathsf{Supp} X = G$ for dimensional reasons.
Let $t$ be a generator of $C_p$, let $x$ and $y$ be non-zero elements of degrees $(t,1)$ and $(1,t)$ respectively.
Again, $x$ and $y$ are invertible and we normalize them in order to have $x^p=y^p=1$. They do not commute, since otherwise the algebra would be the commutative algebra $k(C_p\times C_p)$. In fact $xy$ and $yx$ are both non-zero and they have common degree $(t,t)$. Hence they differ by a scalar: $yx=qxy$. Moreover $q^p=1$ since $x=y^px=q^pxy^p=q^px$. Then $q$ is a primitive root of unity and the grading corresponds to the grading of Proposition \ref{cncn}. \qed

\begin{thm}
Let $k$ be an algebraically closed field of characteristic zero, and let $p$ be a prime. Then
\[ \pi_1 M_p(k) \simeq F_{p-1}\times C_p.\]
\end{thm}
The proof is completely analogous to the proof of Theorem \ref{23}.

We end this section with a computation of the fundamental group of triangular matrix algebras, based on the work of
A. Valenti and M.V. Zaicev \cite{vaza}.

A grading of an upper triangular matrix algebra $T_n(k)$ is \textbf{good} if the elementary matrices ${}_jE_i$ are homogeneous.
Clearly any good grading is completely determined by assigning group elements to sub-diagonal elementary matrices ${}_{i+1}E_i$,
since the idempotents ${}_iE_i$ have necessarily trivial degree.
In other words a good grading is determined as before by a map $m:\{1,\dots,n-1\}\to G$.
The grading is connected if and only if $\Im m$ generates $G$. As before, any good connected grading is a quotient
of the grading given by the free group $F_{n-1}$ on a set $\{s_1, \dots , s_{n-1}\}$ and a map $m$ such that $\Im m = \{s_1, \dots , s_{n-1}\}$.

The main result of \cite{vaza} states that any grading of a triangular algebra is good, without any hypothesis concerning the field
(Theorem 7). As an immediate consequence we obtain the following:

\begin{thm}
Let $k$ be a field and let $T_n(k)$ be the algebra of triangular matrices of size $n$. Then
\[ \pi_1 T_n(k) \simeq F_{n-1}.\]
\end{thm}

\section{\sf Fundamental group of truncated polynomial algebras}

In this section we compute the fundamental group of the group algebra of the cyclic group of
order $p$ in characteristic $p$, in other words we compute the fundamental group of $k[x]/(x^p)$.

\begin{pro}
Let $G$ be a finite group and let $k$ be any field. The usual $G$-grading of the group algebra $kG$
is simply connected.
\end{pro}
\noindent \textbf{Proof:} The Galois covering $kG \# G$ has $G$ as set of objects and, given $s,t \in G$,
\[ {}_t(kG \# G)_s= k(t^{-1}s) . \]
The composition is given by the product of $G$.
In other words, all $k$-vector spaces of morphisms are one-dimensional
and all the structure constants are $1$.
By Lemma \ref{q}, this category is simply connected.\qed

\begin{rem}
As a consequence of the previous proof we recover the Cohen-Mont\-gomery duality theorem for coactions \cite{como}:
the algebras $kG\#G$ and $M_{|G|}(k)$ are isomorphic.
The algebra associated to a finite object category is
obtained as the direct sum of all the vector spaces of morphisms. In particular,
if {all} the vector spaces of morphisms are one-dimensional, we get the matrix algebra.
Hence the algebra associated to the category $kG\#G$  is $M_{|G|}(k)$.
On the other hand, it has been proved in \cite{cima} that the algebra corresponding to the categorical smash product by a
finite group is precisely the usual smash product algebra.
\end{rem}

\medskip

Next we provide an example of a path $k$-algebra of a quiver with admissible relations,
which does not admit a universal cover when the field is of characteristic $p$. The quiver is a loop, and
the relation is given by the $p^{th}$-power of the loop.
There are at least two simply connected coverings by smash categories. One of them is not a covering of
"quivers with relations" in the sense of \cite{ga}.

\begin{pro}
Let $k$ be a field of characteristic $p$. The truncated polynomial algebra $k[x]/(x^p)$ does not admit a universal covering.
\end{pro}
\noindent \textbf{Proof:} First, note that $k[x]/(x^p)$ is isomorphic to the $k$-group algebra of the cyclic group $C_p$ of order $p$,
hence the preceding proposition provides a simply connected covering with group $C_p$.
Note that this covering is the category with $p$ vertices, where all vector spaces of morphisms are one-dimensional
and all the structure constants are $1$. \\
On the other hand, consider the usual $\ZZ$-grading of $k[x]$. Since $(x^p)$ is an homogeneous ideal
- this holds in any characteristic - it induces a grading in $k[x]/(x^p)$.
For this grading, $\left[k[x]/(x^p)\right] \# \ZZ $ is the category which has $\ZZ $ as set of objects,
one-dimensional vector spaces of morphisms from $i$ to $j$ if $0\le j-i < p$, and $0$ otherwise.
In other words the morphisms in the category are generated by morphisms from $i$ to $i+1$ for each integer $i$,
with relations such that any composition of $p$ generators is zero.
As a consequence of this description, each grading of $\left[k[x]/(x^p)\right]\#\ZZ$ is freely determined by assigning a degree to
the one dimensional vector space of morphisms from $i$ to $i+1$.
Hence any homogeneous non-zero closed walk has trivial degree. \\
Recall that by definition of connected grading, any element of the group should be the degree of an homogeneous
walk between objects.
Then the unique group that grades this smash product category  in a connected way is the trivial one.
As a consequence, this covering category is simply connected.
Finally note that the Galois coverings are not isomorphic since their groups of automorphisms are not isomorphic.
In this way we have constructed two non-isomorphic simply connected coverings. \qed

It is well known and easy to prove that the trivial homogeneous component of any grading always contains the ground field $k$.

A grading is called \textbf{fine} if the dimension of each homogeneous component is at most one, see for instance \cite{baseza}.

\begin{thm}
Let $k$ be a field of characteristic $p$ and let $A=k[x]/(x^p)$.
There are two types of connected gradings of $A$,
with no common quotient except the trivial one. The first type corresponds to the group algebra case and the grading group is
$C_p$. In the second one, the grading group is either $\mathbb{Z}$ or any of its quotients.
\end{thm}
\noindent \textbf{Proof:} Let $X$ be a connected basic grading of $A$. There are two cases, according to the existence of an invertible
homogeneous element of nontrivial degree.
First we suppose that there exists an invertible homogeneous element $a$ of degree $s\neq 1$.
We write $ a= a_0 + a_+$ where $a_0\in k^*$ and $a_+\in (x)$, and we normalize $a$ in order to have $a_0=1$.
Since the characteristic of $k$ is $p$, we obtain that $a^p=1$ and $p$ is the order of $a$.
For $i<p$ we infer that $a^i\neq 0$, then $X^{s^i}A\neq 0$.
Moreover $X^{s^i}A\neq X^{s^j}A$ for $i\neq j$, $i,j<p$. Also $1=a^p\in X^{s^p}A$ implies $s^p=1$.
Computing dimensions and since the grading is connected, we deduce that the group is cyclic of order $p$, and the grading is fine. \\
As a second case, assume all homogeneous elements of nontrivial degree belong to the maximal ideal $(x)$, that is
\[\bigoplus_{s\in G, s\neq 1} X^sA\subseteq(x).\]
Consider now the usual valuation $\nu$ on $A$, namely for $f\neq 0$ we have that $\nu(f)$ is the smallest exponent of $x$
appearing in $f$. Of course $\nu(f)=0$ if and only if $f$ is invertible. The valuation $\nu$ has the following properties:
\begin{itemize}
  \item $\nu(f+g) \geq \mathsf{inf}\{\nu(f), \nu(g)\}$ for $f, g, f+g \neq 0.$
  \item $\nu(fg)= \nu(f) + \nu(g)$ for $f,g, fg \neq 0.$
\end{itemize}
Then for $f\neq 0$ we obtain $f=x^{\nu(f)}u$ where $u$ is invertible. \\
Assume first that there exists a homogeneous $g_1 \in X^1A$ of valuation $1$, that is, $g_1=x+u$ with $u \in (x^2)$.
Since $g_1^{p-1}=x^{p-1}$ and $g_1^{p-1}$ is homogeneous, we infer that $x^{p-1}$ is homogenous of degree $1$.
Now, $g_1^{p-2}=x^{p-2}+\lambda x^{p-1}$, so, $x^{p-2}=g_1^{p-2}-\lambda x^{p-1}$ hence $x^{p-2}$ is homogenous of degree $1$.
If we continue with this procedure, we finally get that $x$ is homogenous of degree $1$ and the grading is trivial.
Finally, assume $\nu(g_1)\geq 2$ for any homogeneous $g_1 \in X^1A$.
We claim that there exists a homogeneous $f$ of valuation $1$. If not, for any $g\in(x)$ we have $\nu(g)\geq 2$,
by decomposing $g$ as a sum of its homogeneous components and using the above property of a valuation, which is clearly false since $\nu(x)=1$.
Now $\nu(f^i)=i$ for $i<p$. Since $f^i\in X^{s^i}A$, the latter is not zero.
For dimensional reasons we infer that the support of the grading is $\{1,s,\dots, s^{p-1}\}$ which generates a cyclic group. \qed

\begin{cor}
Let $k$ be a field of characteristic $p$. Then $$\pi_1 \left(k[x]/(x^p)\right) = \mathbb{Z}\times C_p.$$
\end{cor}

\section{\sf Fundamental group of diagonal algebras}

Let $E$ be a finite set and $k$ a field. The diagonal algebra $k^E$ is the vector space of maps from $E$ to $k$ with pointwise
multiplication.
Next we will consider connected gradings of diagonal algebras (see \cite{da} and \cite{bi}). The following result shows
that any abelian group with the  cardinality of a given  set
grades the diagonal algebra in a connected way, provided the field contains enough roots of unity.

\begin{pro}\label{5.11}
Let $E$ be a finite set of order $n$ and let $k$ be a field with enough $n$-th roots of unity. Let $G$ be any abelian group
of order $n$. Then there is a simply connected $G$-grading of $k^E$.
\end{pro}
\noindent \textbf{Proof:} We first sketch the proof of the following well known result.
Let $G$ be any abelian group of order $n$, $E$ a set of cardinal $n$ and $k$ a field containing $n$ different $n$-th roots of unity,
then the algebras $kG$ and $k^E$ are isomorphic. First assume that $G$ is cyclic.
Let $t$ be a generator of $G$ and let $\mu_n$ be the set of $n$-th roots of unity in $k$.
Note that under our assumptions $p$ does not divide $n$, in case $k$ is a field of characteristic $p>0$. Then the set
\[\left\{ e_\zeta = \frac{1}{n}\sum_{i=0}^{n-1}\zeta^it^i\right\}_{\zeta\in\mu_n}\]
is a complete set of orthogonal idempotents of $kG$, which has $n$ elements. This set provides a new basis of $kG$,
proving that $kG$ is isomorphic to $\oplus_{\zeta\in\mu_n}ke_\zeta$, which in turn is identified with $k^E$ through a
bijection between $E$ and $\mu_n$ by considering the Dirac masses in $k^E$. \\
For an arbitrary abelian group $G$ of order $n$, note that $G$ is a direct product of finite cyclic groups.
Note also that a  group algebra $k(G_1\times G_2)$ is isomorphic to $kG_1\otimes kG_2$, while the algebras
$k^{E_1\times E_2}$ and $k^{E_1}\otimes k^{E_2}$ are also isomorphic. The previous case provides the required  isomorphism. \\
Next we prove the statement of the proposition. Consider the algebra $k^E$, an arbitrary abelian group $G$ of order $n$, and an algebra
isomorphism between $kG$ and $k^E$ as before.
The usual $G$-grading of $kG$  provides a grading of $k^E$ by transporting the structure through the isomorphism.
Consequently any abelian group of order $n$  provides a simply connected grading of the algebra $k^E$.\qed

\begin{cor}
Let $n$ be a non-square free positive integer and let $k$ be a field as above. The algebra $k^n$ does not admit a universal covering.
\end{cor}
\noindent \textbf{Proof:} If $n$ is not square free there exist at least two non-isomorphic groups of order $n$.
The result above provides at least
two non-isomorphic  simply connected coverings, then $k^n$ does not admit a universal cover. Moreover,
each abelian group $G$ of order $n$ provides a simply connected grading through the isomorphism of $k^n$ with $kG$.  \qed

The following result is based on the fact that $k\times k$ admits precisely one connected grading. We provide a proof of this, which is also a particular case of S. D\u asc\u alescu's classification in \cite{da} (see also J. Bichon approach in \cite{bi}).

\begin{pro}\label{dos}
Let $k$ be a field of characteristic different from $2$. The fundamental group $\pi_1(k\times k)$ is cyclic of order $2$.
\end{pro}
\noindent \textbf{Proof:} Let $X$ be a connected $G$-grading of $k\times k$ for some group $G$.
The trivial homogeneous component $X^1(k\times k)$ contains the unit of the algebra.
If $X^1(k\times k) = k\times k$, the group is trivial since the support of $X$ is just the trivial element of $G$
and the grading is connected. Otherwise there is exactly one more non-zero homogeneous component $X^s(k\times k)$ which is one dimensional.
Note that $s$ has to generate $G$. We will prove that $s$ is of order $2$. Let $(x,y)$ be a non-zero element of degree $s$.
Clearly $(x,y)^2\neq 0$, moreover $(x,y)^2\in X^{s^2}(k\times k)$. Since there are only two homogeneous components, we infer that
$X^{s^2}(k\times k)=X^1(k\times k)$ or $X^{s^2}(k\times k)=X^{s}(k\times k)$.
In the first case $s^2=1$, while in the second case $s=1$. Consequently there are precisely two connected gradings and the fundamental
group is cyclic of order two.\qed

\begin{lem}\label{productolibre}
Let $A$ and $B$ be algebras provided with connected $G_A$ and $G_B$-gradings $X$ and $Y$ respectively.
Then the algebra $C=A\times B$ has a natural $G_A * G_B$-connected grading $Z$.
As a consequence all quotients of $G_A * G_B$ grade $C$ connectedly.
\end{lem}
\noindent \textbf{Proof:} Consider the following subspaces of $C$:
\[ \begin{array}{ll}
Z^{1}C= X^1A\times Y^1B, &\\
Z^{s}C= X^sA\times 0, \qquad & \mbox{if $s\neq 1$ and $s\in G_A$},\\
Z^{t}C=0\times Y^tB, &\mbox{if $t\neq 1$ and $t\in G_B$},\\
Z^{w}C=0, & \mbox{in the remaining cases.}
\end{array}\]
Note that the support of the grading $Z$ is the union of the supports of $X$ and $Y$. These supports generate $G_A$ and $G_B$ respectively,
hence the support of $Z$ generates $G_A * G_B$.\qed

\begin{ejem}\label{seis}
Let $E_5$ be a set with five elements. There exists a connected $C_6$-grading of $k^{E_5}$.

Indeed let $E_2$ and $E_3$ be sets with two and three elements respectively.
Then $k^{E_5}\cong  k^{E_2}\times k^{E_3}$ and we consider the previous fine and connected gradings given by $C_2$ and $C_3$
of $k^{E_2}$ and $k^{E_3}$ respectively. The previous lemma shows that the product group $C_2 * C_3$  grades the
product algebra $k^{E_5}$ in a connected way, as well as any of its quotients, in particular $C_6$.
\end{ejem}

This example is in fact the basis of the general procedure developed by S. D\u asc\u alescu in order to describe all the connected gradings of a diagonal algebra. First we rephrase one of his results.

\begin{lem}\cite[Lemma 1]{da}\label{ergodic}
Let $k^n$ be a diagonal algebra. Any connected $G$-grading with one-dimensional trivial homogeneous component is given by the usual $G$-grading of $kG$, where $G$ is any abelian group of order $n$.
\end{lem}

Note that D\u asc\u alescu call \textbf{ergodic} a grading with one-dimensional trivial homogeneous component.
For $n=2$,  a nontrivial grading has to be ergodic, hence we recover the fact that there is only one
nontrivial grading of $k\times k$ as in Proposition \ref{dos}.

D\u asc\u alescu provides a description of all the gradings of $k^E$, which is based on ergodic ones (see \cite[Theorem 5]{da}).
We shall use it in order to compute $\pi_1(k^n)$ for small values of $n$. In order to state his result, we fist consider the following specific connected gradings of a diagonal algebra which are modeled on Example \ref{seis}

Roughly speaking the specific gradings are free product gradings of connected ergodic ones based on a product algebra decompositions of a diagonal algebra. Note that connected ergodic gradings of diagonal algebras are classified by Lemma \ref{ergodic}.

More precisely, let $A=k^E$ be a diagonal algebra and let $M_1,\dots, M_s$ be a partition of $E$. Let $A_{M_i}$ be the algebra $A(e_{M_i})$, where $e_{M_i}=\sum_{x\in M_i}\delta_x$, and $\delta_x$ is the Dirac mass at $x$. It is well known and easy to prove that any direct product decomposition of $k^E$ is obtained in this way. Let $H_i$ be an abelian group of order $\#(M_i)$ and finally let $X_i$ be the corresponding $H_i$-ergodic grading of $A_{M_i}$. Then by Lemma \ref{productolibre}, the group $H_1 * \dots * H_s$ provides a connected grading of $A=A_{M_1}\times \dots \times A_{M_s}$, which we call specific.

\begin{thm}\cite{da}\label{diagonal}
Let $E$ be a finite set and let $k$ be a field containing all roots of unity of order less than or equal to $\# E$.
Any connected grading of $k^E$ is a quotient of a specific grading.
\end{thm}

\begin{cor}
Let $k$ be a field containing all roots of unity of order $2$ and $3$.
Then $\pi_1(k^3)=C_2\times C_3$.
\end{cor}
\noindent \textbf{Proof:} The two nontrivial partitions of $\{1,2,3 \}$ provide connected gradings by $C_2$ and $C_3$.
Clearly they do not have nontrivial common quotients.

\begin{thm}
Let $k$ be a field containing all roots of unity of order $2$, $3$ and $4$.
Then $\pi_1(k^4)=(C_2 * C_2)\times C_4 \times C_2\times C_2 \times C_3 = (C_2 * C_2)\times C_6 \times C_4\times C_2.$
\end{thm}
\noindent \textbf{Proof:} The specific gradings of $k^4$ are given by the partitions of the set $\{1,2,3,4 \}$ as follows.
\begin{center}
\begin{tabular}{|c|c|c|}
  \hline

Group & Dimension   & Dimension \\
 &  of the trivial component & of other components \\

  \hline

 $\{ 1\}$ & 4 & 0\\

 \hline

 $C_2 * C_2$ & $2$ & ${1,1}$\\

 \hline

 $C_3$ & $2$ & $1,1$\\

 \hline

$C_2$ & $3$ & $1$\\

 \hline

$C_4$ & $1$ & $1,1,1$\\

 \hline

$C_2 \times C_2$ & $1$ & $1,1,1$\\

 \hline

\end{tabular}  \end{center}
An inspection of the possible common quotients taking into account the structure of the groups and the dimension of the trivial homogeneous
components shows that the $C_2$-grading is a quotient of the $C_2 * C_2$-grading.
Moreover, there is no other nontrivial common quotient.
\qed


\footnotesize
\noindent C.C.:
\\Institut de Math\'{e}matiques et de Mod\'{e}lisation de Montpellier I3M,\\
UMR 5149\\
Universit\'{e}  Montpellier 2,
\\F-34095 Montpellier cedex 5,
France.\\
{\tt Claude.Cibils@math.univ-montp2.fr}

\noindent M.J.R.:
\\Departamento de Matem\'atica,
Universidad Nacional del Sur,\\Av. Alem 1253\\8000 Bah\'\i a Blanca,
Argentina.\\ {\tt mredondo@uns.edu.ar}

\noindent A.S.:
\\Departamento de Matem\'atica,
 Facultad de Ciencias Exactas y Naturales,\\
 Universidad de Buenos Aires
\\Ciudad Universitaria, Pabell\'on 1\\
1428, Buenos Aires, Argentina. \\{\tt asolotar@dm.uba.ar}


\begin{thebibliography}{99}


\bibitem{al} Aljadeff, E.; Haile, D., Natapov, M. Graded identities of matrix algebras and the universal graded algebra.
To appear in Trans. Am. Math. Soc.
\texttt{arXiv:0710.5568}


\bibitem{baseza} Y. Bahturin, S. Sehgal, M. Zaicev,
Group gradings on associative algebras,
J. Algebra 241 (2001), no. 2, 677--698.

\bibitem{bashe} Y. Bahturin, I. Shestakov,
Gradings of simple Jordan algebras and their relation to the gradings of simple associative algebras,
Comm. Algebra  29 (2001), 4095--4102.

\bibitem{baza} Y. Bahturin, M. Zaicev,
Group gradings on matrix algebras,
Canad. Math. Bull. 45 (2002), no. 4, 499--508.

\bibitem{bi} J. Bichon, Algebraic quantum permutation groups,
 Asian-Eur. J. Math.  1  (2008),  1--13.

\bibitem{boboc} C. Boboc,
Gradings of matrix algebras by the Klein group,
Comm. Algebra  31 (2003), 2311--2326.

\bibitem{bodas} C. Boboc, S. D\u asc\u alescu,
Gradings of matrix algebras by cyclic groups,
Comm. Algebra 29 (2001), 5013--5021.

\bibitem{bodas06} C. Boboc, S. D\u asc\u alescu,
Good gradings of matrix algebras by finite abelian group of prime index,
Bull. Math. Soc. Sc. Math. Roumanie 49 (2006), 5-–11.

\bibitem{bodas07} C. Boboc, S. D\u asc\u alescu,
Group gradings on $M_3(k)$,
Comm. Algebra  35 (2007), 2654–-2670.

\bibitem{boga} K. Bongartz, P. Gabriel,
Covering spaces in representation-theory,
Invent. Math. 65 (1981/82), 331--378.

\bibitem{caedasnas} S. Caenepeel, S. D\u asc\u alescu, C. N\u ast\u asescu,
On gradings of matrix algebras and descent theory,
Comm. Algebra 30 (2002), 5901--5920.

\bibitem{cima} C. Cibils, E. Marcos,
Skew category, Galois covering and smash product of a $k$-category,
Proc. Amer. Math. Soc. 134 (2006),  no. 1, 39--50.

\bibitem{cireso} C. Cibils, M. J. Redondo, A. Solotar,
The intrinsic fundamental group of a linear category,
\texttt{arXiv:0706.2491}

\bibitem{chunlee} J. Chun, J. Lee,
Abelian group gradings on full matrix rings,
Comm. Algebra  35 (2007), 3095--3102.

\bibitem{como}M. Cohen, S. Montgomery,
Group-graded rings, smash products, and group actions,
Trans. Amer. Math. Soc. 282  (1984), 237--258.

\bibitem{da} S. D\u asc\u alescu,
Group gradings on diagonal algebras,
Arch. Math. (Basel) 91 (2008), 212--217.

\bibitem{dinr} S. D\u asc\u alescu, B. Ion, C. N\u ast\u asescu, J. R\'\i os Montes,
Group gradings on full matrix rings,
J. Algebra 220 (1999), no. 2, 709--728.

\bibitem{ga} P. Gabriel,
The universal cover of a representation-finite algebra,
Representations of algebras (Puebla, 1980), pp. 68--105,
Lecture Notes in Math., 903, Springer, Berlin-New York, 1981.

\bibitem{green} E. Green,
Graphs with relations, coverings and group-graded algebras,
Trans. Amer. Math. Soc.  279  (1983),  no. 1, 297--310.

\bibitem{grma} E. Green, E. Marcos,
Graded quotients of path algebras: A local theory,
J. Pure Appl. Alg. 93 (1994), 195–226.

\bibitem{kabodas} R. Khazal, C. Boboc, S. D\u asc\u alescu,
Group gradings of $M\sb 2(K)$,
Bull. Austral. Math. Soc. 68 (2003), no. 2, 285--293.

\bibitem{le} P. Le Meur,
The universal cover of an algebra without double bypass,
J. Algebra  312  (2007),  no. 1, 330--353.

\bibitem{mi} B. Mitchell,
Rings with several objects,
Adv. Math.  8 (1972), 1--161.

\bibitem{vaza} A. Valenti, M. Zaicev,
Group gradings on upper triangular matrices,
Arch. Math. (Basel)  89  (2007), 33--40.
		
\end{thebibliography}
\end{document}